\documentclass[reqno,11pt]{amsart}
\usepackage{amscd,amssymb}
\usepackage{latexsym}
\usepackage[all]{xy}

\newtheorem{theorem}{Theorem}[section]
\newtheorem{proposition}[theorem]{Proposition}
\newtheorem{lemma}[theorem]{Lemma}

\theoremstyle{definition}

\newtheorem{remark}[theorem]{Remark}

\def\into{\hookrightarrow}
\def\eps{\epsilon}
\def\del{\delta}
\def\alp{\alpha}
\def\Sp{{\rm Sp}}
\def\bfR{{\mathbf R}}
\def\Quot{{\rm Quot}}
\def\calF{{\mathcal F}}
\def\bbN{{\mathbb N}}
\def\bbP{{\mathbb P}}
\def\frm{{\mathfrak m}}
\def\calOcirc{{\mathcal O^\circ}}

\begin{document}

\title{Power series over generalized Krull domains}
\author{Elad Paran}
\address{School of Mathematics\\
Tel Aviv University\\
Ramat Aviv, Tel Aviv 69978, Israel}
\email{paranela@post.tau.ac.il}
\author{Michael Temkin}
\address{Department of Mathematics\\
University of Pennsylvania\\
Philadelphia, PA  19104 USA} \email{temkin@math.upenn.edu}

\begin{abstract}
We resolve an open problem in commutative algebra and Field
Arithmetic, posed by Jarden -- Let $R$ be a generalized Krull
domain. Is the ring $R[[X]]$ of formal power series over $R$ a
generalized Krull domain? We show that the answer is negative.
\end{abstract}
\maketitle
\section{Introduction}
Recall \cite[\S13.VI.13]{ZS} that an integral domain $R$ is called
a {\bf Krull domain}, if there exists a family $\calF$ of {\bf
discrete} rank-1 valuations of $K = \Quot(R)$, satisfying the
following properties:

(a) For each $v \in \calF$, the valuation ring $R_v$ of $v$ in $K$ is the localization of $R$ with respect to
$\frm_v = \{a \in R| v(a) > 0\}$.

(b) The intersection of all valuation rings $\cap_{v \in \calF} R_v$ is $R$.

(c) For each $0 \neq a \in R$, $v(x) = 0$ for all but finitely many $v \in \calF$.

Every Noetherian integrally closed domain is a Krull domain
\cite[Theorem 12.4(i)]{Matnew}.

Krull domains play important role in commutative algebra. For
example, it is known that the integral closure of a Noetherian
domain is (unfortunately) not necessarily Noetherian. However, by
the Mori-Nagata integral closure theorem it is necessarily a Krull
domain \cite[{\S}A.41]{Matold}.

If $R$ is a Krull domain, so is the ring of polynomials $R[X]$, as
well as the ring of formal power series $R[[X]]$ \cite[Theorem
12.4(iii)]{Matnew}.

In 1981 Weissauer \cite[\S7]{Ws} introduced the notion of a {\bf
generalized Krull domain} -- a domain $R$ is called a generalized
Krull domain \cite[\S15.4]{FJ}, if it is equipped with a family
$\calF$ of {\bf real (= rank-1)} valuations (not necessarily
discrete), satisfying the same conditions (a),(b),(c) above. The
importance of generalized Krull domains in Field Arithmetic and
Galois theory lies in Weissauer's theorem -- the quotient field of
a generalized Krull domain of dimension exceeding 1 is Hilbertian
\cite[Theorem 15.4.6]{FJ}. This widely general theorem provides
many non-trivial Hilbertian fields, and has had extensive use in
recent results in Field Arithmetic, concerning Galois theory over
quotient fields of complete domains (e.g. \cite{Par}).

If $R$ is a generalized Krull domain, then so is $R[X]$. However, up until now it was unknown \cite[Problem
15.5.9(a)]{FJ} if the same holds for $R[[X]]$, as in the discrete case.

In this paper we prove that the answer is negative. We show that
if $R$ is a {\bf complete} non-Noetherian real valuation ring (in
particular, $R$ is a generalized Krull domain), then $R[[X]]$ is
never a generalized Krull domain.

The ring $R[[X]]$ exhibits certain weird behaviors which makes it
difficult to analyze directly. For example, although the Krull
dimension of $R$ is 1, the Krull dimension of $R[[X]]$ is not 2
(as it is in the discrete case), but infinite.

To overcome the wild behavior of this ring, we embed it into 
the ring $R\{X\}$ of convergent power series over $R$, by using substitutions of the form $X \mapsto a X$, for
$a \in R$ with positive valuation. The latter ring is well known from rigid analytic geometry, and has pleasant
properties that we exploit (e.g. the Weierstrass preparation theorem holds in this ring). Intuitively, by
choosing elements $a \in R$ with valuation tending to 0, we approximate $R[[X]]$ better and better by copies of
$R\{X\}$, which enable us to use the properties of $R\{X\}$ to gain information on $R[[X]]$. This interplay
allows us to prove the main result of this paper -- Theorem \ref{mainth}.


A consequence of our result is that one cannot apply Weissauer's
theorem to prove that $\Quot(R[[X]])$ is Hilbertian. This remains
an open question \cite[Problem 15.5.9(b)]{FJ}, and we hope this
work is a step towards its resolution.

\section{Power series over a complete real valuation ring}

Let $K$ be a field equipped with a non-archimidean real valuation
$w$, which is not discrete, and let $R$ be the valuation ring of
$K$. Equivalently, $R$ is a non-Noetherian integrally closed local
domain of dimension 1. In particular, $R$ is a generalized Krull
domain (where the corresponding family of valuations is just
$\{w\}$). Put $D = R[[X]]$. Let $\frm$ be the maximal ideal of
$R$. We assume $K$ is {\bf complete} with respect to $w$. We rely
on simple algebraic properties of the ring $K\{X\}$ of convergent
power series over $K$ (equivalently, this is the ring of
holomorphic functions on the unit disc in $\bbP ^1 _K$). These are
developed in a short and self-contained manner in \cite[\S1]{HV}.

In the following remark we give a rigid-geometric interpretation of the ring $R[[X]]$, and explain how one can
use rigid geometry to establish its properties. This will not be used in the sequel and can be skipped by a
reader with low rigid-geometric motivation.

\begin{remark}
In this remark only we use a (multiplicative) absolute value $|\
|:K\to\bfR_+$ instead of the (additive) valuation $w$. This is the
common choice in non-Archimedean geometry. We have obvious
inclusions $R\{X\}\into R[[X]]\into R\{\frac{X}{a}\}$, where
$a\in\frm$. On each $K\{\frac{X}{a}\}$ we have the Gauss valuation
$|\ |_a$ given by $|\sum_{i=0}^\infty f_iX^i|_a= \max _{i \in
\bbN}|a|^i|f_i|$. Consider the unit rigid disc $B=\Sp(K\{X\})$
with center at zero and of radius one. For each $a\in R$ we
consider a smaller disc $B_a=\Sp(K\{\frac{X}{a}\})$ given by
$|X|\le a$, then
$$|f(X)|_a=\max_{c\in B_a}|f(c)|=\max_{\alp\in
K^a,|\alp|\le|a|}|f(\alp)|$$ by \cite[5.1.4/6]{BGR} (i.e. the
Gauss norm is the supremum norm on the disc). The disc
$B^-:=\cup_{a\in\frm} B_a$ is a non-quasi-compact subdomain of $B$
given by $|X|<1$. If $\calOcirc$ denotes the sheaf of functions of
absolute value $\le 1$, then $\calOcirc(B_a)=R\{\frac{X}{a}\}$ and
therefore $\calOcirc(B^-)=\cap_{a\in\frm}R\{\frac{X}{a}\}=R[[X]]$,
i.e. $R[[X]]$ is the ring of functions whose norm does not exceed
$1$ on $B^-$.

(i) Extend $|\ |$ from $R\{X\}$ to $R[[X]]$ by setting $|f|^-:=\sup_i |f_i|$ for any $f=\sum_i f_iX^i\in
R[[X]]$. Obviously, $|f|^-=\sup_{a\in\frm}|f|_a=\lim_{|a|\nearrow 1}|f|_a$. Since $|f|_a$ are multiplicative,
$|f|^-$ is multiplicative and hence is a valuation. Note also that it follows that $|f|^-=\sup_{c\in
B^-}|f(c)|$.

(ii) Next we claim that if $f(c)=0$ for $c\in\frm$, then $f(X)/(X-c)\in R[[X]]$. If $c=0$ then this is clear,
and the general case reduces to this one by the coordinate change $X'=X-c$ (any point of a non-Archimedean disc
is its center).

(iii) Furthermore, any point $c\in B^-$ corresponds to a uniquely
defined irreducible monic polynomial $g(X)\in K[X]$ whose roots in
$K^a$ are of absolute value strictly smaller than $1$. The latter
happens iff $g(X)\in X^n+\frm[X]$. We claim that if $f(X)\in
K[[X]]$ vanishes at $c$, then $\frac{f(X)}{g(X)}\in K[[X]]$.
Indeed, we proved this in (ii) for linear polynomials and the
general case is obtained by embedding $K[[X]]$ into $K^a[[X]]$
(and some additional care for the inseparable case).

(iv) Any function $f(X)\in K\{\frac{X}{a}\}$ has $n$ geometric zeros on $B_a$, where $n:=\deg_{B_a}(f)$ is the
maximal integer so that $|f(X)|_a=|f_nX^n|$. (For example, use the Weierstrass division theorem
\cite[5.2.2/1]{BGR}).

(v) If $f(X)\in\frm[[X]]$ satisfies $|f|^-=1$, then
$\lim_{|a|\nearrow 1}\deg_{B_a}(f(X))=\infty$. Using (iii) and
(iv) we obtain that such $f(X)$ is divided by polynomials
$g_1(X)$, $g_1(X)g_2(X)$, $g_1(X)g_2(X)g_3(X)$, etc., where the
$g_i$'s are as in (iii) and if $\alp_i$ denotes a root of $g_i$
then $|\alp_i|\nearrow 1$.

(vi) The phenomenon from (v) is possible because $B^-$ is not
quasi-compact, so a function $f(X)$ on $B^-$ can have infinitely
many zeros though it has finitely many zeros on each quasi-compact
subdomain $B_a$. \qed \end{remark}

We now give a direct proof of the properties of $D = R[[X]]$
described in the remark, and then we use these properties to show
that $D$ cannot be a generalized Krull domain.

\begin{lemma}\label{extlem}
The valuation $w$ extends to a valuation on $D$ given by the formula $w(\sum f_i X^i)) = \inf w(f_i)$. 
\end{lemma}
\begin{proof}
The only non-trivial part is to check that $w((f \cdot g) (X)) = w(f(X)) + w(g(X))$. Note that
$w(f(X)) = \lim_{\del \to 0^+} (\min w(f_i) + i\del)$. Indeed, the left hand side is clearly not greater than
the right hand side. Conversely, let $\eps > 0$, and choose $n \in \bbN$ such that for each $i \geq n$ we have
$w(f_i) - w(f(X)) < \eps$. In particular ${\eps + w(f(X)) - w(f_n) \over n} > 0$, so we may choose $0 < \del$
with $w(f_n) + n\del < \eps + w(f(X))$, hence $\min (w(f_i) + i \del) < \eps + w(f(X))$.

Now, for each $a \in \frm$ we have $f(aX),g(aX) \in K\{X\}$, where
$K\{X\} = \{\sum a_i X^i \in K[[X]]| w(a_i) \to \infty\}$
\cite[\S1]{HV}. Hence $w((fg)(aX)) = w(f(aX)) + w(g(aX))$, by
\cite[Lemma 1.3(i)]{HV}. Then

$$w(fg(X)) = \lim_{w(a) \to 0^+} w(fg(aX)) =
\lim_{w(a) \to 0^+} w(f(aX)) + w(g(aX)) = $$

$$ \lim_{w(a) \to 0^+} w(f(aX)) + \lim_{w(a) \to 0^+} w(g(aX)) =
w(f(X)) + w(g(X))$$
\end{proof}

Extend $w$ naturally to $F = \Quot(D)$.

\begin{lemma}\label{divlem}
Suppose $c \in \frm, f(X) \in D$, such that $c$ is a root of $f(X)$ of order at least $n$. That is, $f(c) =
f'(c) = \ldots = f^{(n-1)}(c) = 0$. Then $f(X)$ is divisible by $(X - c)^n$ in $D$. Moreover, if $w(f(X)) = 0$
and $f_0,f_1,\ldots $ are of positive valuation, then ${f(X) \over (X - c)^n}$ satisfies these properties as
well.
\end{lemma}
\begin{proof}

If $c = 0$ the claim is obvious. Assume $c \neq 0$ and let $g_i =
-({f_0 \over c^{i + 1}} + {f_1 \over c^i} + \ldots + {f_i \over
c})$ for each $i \geq 0$. Then $g(X) = \sum g_i X^i \in K[[X]]$
satisfies $g(X)(X-c) = f(X)$.

Let $i \geq 1$. Since $f(c)$ converges, we have $w(f_j c^j) \to \infty$, hence $w(f_j c^{j-i}) \to \infty$.
Hence $f_i + f_{i+1} c + f_{i+2}c^2 + \ldots$ converges to an element $a_i \in R$. Then $b_i = a_i c^i = f_i c^i
+ f_{i+1} c^{i+1} + \ldots$ has valuation $w(b_i) \geq iw(c)$. Since $f(c) = 0$, we have $f_0 + f_1 c + \ldots +
f_{i-1} c^{i-1} = -b_i$, hence $w(f_0 + f_1 c + \ldots + f_{i-1} c^{i-1}) \geq iw(c)$. Thus $w(g_{i-1}) \geq 0$
and hence $g(X) \in R[[X]]$.

Now suppose $f_i \in \frm$ for each $i \geq 0$, and $w(f(X)) = 0$. Let $i \geq 1$. Then $w(f_i + f_{i+1} c +
f_{i+2}c^2 + \ldots)
> 0$, hence $w(f_0 + f_1 c + \ldots + f_{i-1} c^{i-1}) > iw(c)$,
so $w(g_{i-1}) > 0$. Moreover, $w(g(X)) = w(f(X)) - w(X - c) = 0 - 0 = 0$.

We have proven the claim for the case $n = 1$. The general case follows by induction.
\end{proof}

Note that the proof of Lemma \ref{divlem} does not rely on the
fact that $K$ is complete.

\begin{remark}\label{reminv} If an element $f(X) = \sum_{i=0}^\infty f_i X^i \in K\{X\}$ is
invertible then $w(f_0) < w(f_i)$ for each $i > 0$.
\end{remark}\begin{proof} Without loss of generality, $f_0 \neq 0$. Let $g(X)$ be the
inverse of $f(X)$ in $K[[X]]$. For each $0 \neq h(X) \in K\{X\}$,
let $\textrm{p.deg}(h(X)) = \max (n \mid w(h_n)=w(h(X)))$ be the
pseudo degree of $h(X)$ \cite[Definition 1.4]{HV}. It is standard
to check that $\textrm{p.deg}(((h_1\cdot h_2)(X))) =
\textrm{p.deg}(h_1(X)) + \textrm{p.deg}(h_2(X))$. In particular,
if $f(X)$ is invertible, then $g(X) \in K\{X\}$, hence $0 =
\textrm{p.deg}((fg)(X)) = \textrm{p.deg}(f(X)) +
\textrm{p.deg}(g(X))$, so $\textrm{p.deg}(f(X)) =
\textrm{p.deg}(g(X)) = 0$.
\end{proof}

The converse of remark \ref{reminv} also holds, but we shall not
need it.

\begin{lemma}\label{poslem}
Suppose $f(X) \in D$ satisfies $w(f(X)) = 0$, and $f_0,f_1,\ldots$ are of positive valuation.
Then there exists a monic irreducible polynomial $q(x) \in R[X]$ of positive degree, such that $q(x)$ divides
$f(X)$ in $R[[X]]$, $w(q_0) > 0$, and such that $g(X) = {f(X) \over q(X)}$ satisfies $w(g(X)) = 0$ and
$g_0,g_1,\ldots$ are of positive valuation.
\end{lemma}
\begin{proof}
Write $f(X) = \sum f_i X^i$. Since $w(f_0) > w(f(X)) = 0$, there
exists $i >0$ such that $w(f_i) < w(f_0)$. Choose $a \in \frm$
with $w(a)$ sufficiently small, such that $w(f_i) + iw(a) <
w(f_0)$. Let $g(X) = f(a X) \in R\{X\}$. Since $w(g_i) < w(g_0)$
the element $g(X)$ is not invertible in $K\{X\}$, by Remark
\ref{reminv}. By \cite[Lemma 1.9]{HV} we may write $g(X) =
r(X)u(X)$, where $u(X)$ is invertible in $K\{X\}$ and $r(X) \in
K[X]$. By multiplying $r(X)$ with an element of $K^\times$ (and
dividing $u(X)$ by it), we may assume $r(X) \in R[X]$. Since
$g(X)$ is not invertible, $r(X)$ must be of positive degree. Let
$p(X) \in R[X]$ be a monic irreducible factor of $r(X)$.

Let $\tilde{K}$ be the algebraic closure of $K$, and extend $w$ to $\tilde{K}$. Let $\tilde{R}$ be the valuation
ring in $\tilde{K}$ lying over $R$, and put $\tilde{D} = \tilde{R}[[X]]$. Extend $w$ further to $\tilde{D}$ by
$w(\sum a_i X^i) = \inf (w(a_i))$.

Let $c_1,\ldots,c_n$ be the distinct roots of $p(X)$ in $\tilde{K}$, and let $e_i$ be the multiplicity of $c_i$,
for each $1 \leq i \leq n$. Then $a c_i$ is a root of $f(X)$ of multiplicity at least $e_i$. It follows by Lemma
\ref{divlem} that $f(X)$ is divisible in $\tilde{R}[[X]]$ by $(X - ac_1)^{e_1} \cdot \ldots \cdot (X -
ac_n)^{e_n} = a^n p({X \over a})$. Note that $q(X) = a^n p({X \over a})$ is a monic polynomial in $R[X]$, and is
irreducible since $p(X)$ is irreducible. Moreover, $w(q_0) \geq nw(a) > 0$.

Write $h(X) = {f(X) \over q(X)}$. By applying Lemma \ref{divlem} consecutively to the elements 
$c_1, c_2, \ldots,c_n$ (in the ring $\tilde{R}[[X]]$), we get $w(h(X)) = 0$ and $w(h_i) > 0$ for each $i \geq
0$. Finally, $h(X) \in \tilde{R}[[X]] \cap K((X)) = R[[X]]$ (where the intersection is taken inside
$\tilde{K}((X))$).
\end{proof}

\begin{proposition}\label{infprop}
Let $f(X) \in D$ such that $w(f(X)) = 0, f_0,f_1,\ldots \in \frm$. Then $f$ has infinitely many factors which
are monic irreducible elements of $R[X]$ with constant term in $\frm$.
\end{proposition}
\begin{proof}
First note that if $p(X) \in R[X]$ is a 
factor of $f(X)$ with $w(p_0) > 0$, then $p(X)$ divides $f(X)$ finitely many times. Indeed, if $p(X)^n | f(X)$
for each $n \geq 1$, then for a fixed element $a \in \frm$ we have $p(aX) | f(aX)$ for all $n \geq 1$. But
$w(p(aX)) > 0$ (since $w(p_0) > 0$), hence $w(f(aX)) = \infty$. Thus $f(a X) = 0$, hence $f(X) = 0$, a
contradiction.

By Lemma \ref{poslem} $f(X)$ has a monic irreducible factor
$p_1(X)$ (in $R[X]$), with constant term in $\frm$. Suppose by
induction that we have constructed $n$ such factors
$p_1(X),\ldots,p_n(X)$. Dividing $f(X)$ by sufficiently large
powers of $p_1(X),\ldots,p_n(X)$ we get an element $g(X) \in D$
such that $p_1(X),\ldots,p_n(X)$ are not factors of $g(X)$.
Moreover, by Lemma \ref{divlem} $w(g(X)) = 0$ and $g_i \in \frm$
for each $i \geq 0$. By Lemma \ref{poslem} $g(X)$ (and hence
$f(X)$) has a monic irreducible factor $p_{n+1}(X) \in R[X]$ with
constant term in $\frm$.
\end{proof}

\begin{lemma}\label{finlem}
Suppose $D$ is a generalized Krull domain, and let ${\calF}$ be a corresponding family of valuations. Let $t \in
\frm$. Then only finitely many valuations in ${\calF}$ are non-trivial on $R$. Moreover, if $v \in {\calF}$ is
non-trivial on $R$, then $v(t) > 0$.
\end{lemma}
\begin{proof}
Suppose $v \in {\calF}$ is non-trivial on $R$, and let $a \in R$
with $v(a) > 0$. For a sufficiently large $n \in \bbN$ we have
$w(t^n) \geq w(a)$, hence $b = {t^n \over a} \in R$. If $v(t) = 0$
then $v(b) < 0$, a contradiction.

Thus $v(t) > 0$ for each $v \in {\calF}$ which is non-trivial on $R$. Hence only finitely many such valuations
exist in ${\calF}$.
\end{proof}

\begin{theorem}\label{mainth}
The ring $D = R[[X]]$ is not a generalized Krull domain.
\end{theorem}
\begin{proof}
Suppose $D$ is a generalized Krull domain, and let ${\calF}$ be a corresponding family of valuations. Let
${{\calF}}_0$ be the subfamily of valuations in ${\calF}$ which are non-trivial on $R$. By Lemma \ref{finlem}
${\calF}_0$ is finite.

Let $p(X)$ be a monic irreducible element of $R[X]$ with constant
term in $\frm$. Then there exists $v_p \in {\calF} \setminus
{\calF}_0$ with $v_p(p(X))
> 0$. Indeed, suppose $v(p(X)) = 0$ for each $v \in {\calF}
\setminus {\calF}_0$. Fix $t \in \frm$. For each $v \in {\calF}_0$
we have $v(t) > 0$, by Lemma \ref{finlem}. Since ${\calF}_0$ is
finite, for a sufficiently large $n \in \bbN$ we have $v(t^n) >
v(p(X))$, for each $v \in {\calF}_0$. Put $a = t^n, h(X) = {a
\over p(X)} \in F$. Then $v(h(X)) > 0$ for all $v \in {\calF}_0$
and $v(h(X)) = 0$ for all $v \in{\calF} \setminus {\calF}_0$. Thus
$h(X) \in \bigcap_{v \in {\calF}} D_v$, hence $h(X) \in D$, by our
assumption. Since $p(X)$ is monic, $w(p(X)) = 0$, hence $w(h(X)) =
w(a)$. It follows that $w(h_i) \geq w(a)$ for each $i \geq 0$, so
${1 \over p(X)} = {h(X) \over a} \in R[[X]]$. Thus we must have
$w(p_0) = 0$, a contradiction. This proves the existence of $v_p$.

Next, note that for each two distinct such polynomials $p, q$, we
have $v_p(q) = 0$. Indeed, since $p, q$ are irreducible, they are
co-prime, so there exist $r(X),s(X) \in K[X]$ such that $p(X)r(X)
+ q(X)s(X) = 1$. By multiplying with a non-zero element of $R$, we
may assume that $r(X),s(X) \in R[X]$ satisfy $p(X)r(X) + q(X)s(X)
\in R \setminus \{0\}$. Since $v_p$ is trivial on $R$, we have
$v_p(p(X)r(X) + q(X)s(X)) = 0$. Since $v_p(p(X)) > 0$, we must
have $v_p(q(X)) = v_p(s(X)) = 0$. It follows that if $p, q$ are
distinct, so are $v_p,v_q$.

Finally, choose an element $f(X) \in D$ with $w(f(X)) = 0$ and
\newline \noindent $f_0,f_1,\ldots \in \frm$. By Proposition \ref{infprop} $f(X)$
has infinitely many monic irreducible factors
$p_1(X),p_2(X),\ldots \in R[X]$ with constant term in $\frm$. Then
$v_{p_i}(f(X)) \geq v_{p_i}(p_i(X)) > 0$ for each $i \geq 1$. Thus
$v(f(X)) > 0$ for infinitely many $v \in {\calF}$, a
contradiction.
\end{proof}


\begin{thebibliography}{BGR86}
\bibitem[BGR]{BGR}
Bosch, S.; G\"untzer, U.; Remmert, R.: {\it Non-Archimedean analysis. A systematic approach to rigid analytic
geometry}, Springer, Berlin-Heidelberg-New York, 1984.

\bibitem[FrJ]{FJ}
Fried, M.D.; Jarden, M.: {\it Field Arithmetic, 2nd edition,
revised and enlarged by M. Jarden}, Ergebnisse der Mathematik III
{\bf 11}, Sprin\-ger 2005.

\bibitem[FrP]{FP}
Fresnel, J.; van der Put, M.: {\it  Rigid Analytic Geometry and Its Applications}, Progress in Mathematics {\bf
218}, Birkh\"auser Boston, 2004.

\bibitem[HaV]{HV}
Haran, D.; V\"olklein, H.: {\it Galois groups over complete valued fields}, Israel Journal of Mathematics {\bf
93} (1996), 9--27.

\bibitem[Mat1]{Matold}
Matsumura, H.: {\it Commutative Algebra}, second edition, The Benjamin/Cummings Publishing Company, 1980.

\bibitem[Mat2]{Matnew}
Matsumura, H.: {\it Commutative ring theory}, Cambridge University Press, 1986.

\bibitem[Par]{Par}
Paran, E.: {\it Split embedding problems over complete domains},
Annals of mathematics, to appear.

\bibitem[Ws]{Ws}
Weissauer, R.: {\it Der Hilbertsche Irreduzibilit\"atssatz},
Journal f\"ur die reine und angewandle Mathematik {\bf 334}
(1982), 203--220.

\bibitem[ZaS]{ZS}
Zariski, O.; Samuel, P.: {\it Commutative Algebra}, Vol. II, \rm van Nostrand 1960.

\end{thebibliography}
\end{document}